\documentclass[11pt]{article}
\usepackage{amssymb}

\newtheorem{theorem}{Theorem}
\newtheorem{lemma}{Lemma}

\title{$k$-Term Arithmetic Progressions in Very Sparse Sumsets}
\author{Ernie Croot}

\begin{document}

\maketitle

\section{Introduction}

One of the main focuses in combinatorial (and additive) number theory
is that of ``understanding'' the structure of the sumset 
$A+B = \{a+b : a \in A, b \in B\}$, given certain 
information about the sets $A$ and $B$.  For example, one such problem is
to determine the length of the longest arithmetic progression in this 
sumset, given that $A,B \subseteq \{0,1,2,...,N\}$ and $|A|,|B| > \delta N$,
for some $0 < \delta \leq 1$.  
The first major progress on this problem was due to J. Bourgain 
\cite{bourgain}, who proved the beautiful result:

\begin{theorem}  If $A,B \subseteq \{0,1,...,N\}$ and $|A| = \gamma N$
and $|B| = \delta N$, then for $N$ large enough, the set 
$A+B$ contains an arithmetic progression of length 
$L  > \exp[c (\gamma \delta \log N)^{1/3} - \log\log N]$, for some 
constant $c$.  
\end{theorem}

Then, I. Ruzsa \cite{ruzsa} gave an ingenious construction, which is
the following theorem:

\begin{theorem} For every $\epsilon > 0$ and every sufficiently large 
prime $p$, there exists a symmetric set $A$ of residues modulo $p$
with $|A| \geq p(1/2 - \epsilon)$, such that $A+A$ contains no arithmetic
progression modulo $p$ having length $\geq \exp( (\log p)^{2/3+\epsilon})$.
\end{theorem}

A simple consequence of this theorem is that for $N$ sufficiently large,
there exists a subset $A$ of the integers in $[1,N]$ with 
$|A| > (1/4 - \epsilon) N$, having no arithemtic progressions of length
$\geq \exp( (\log N)^{2/3+\epsilon})$, which shows that the $1/3$ in
Bourgain's result cannot be improved to any number beyond $2/3$.  

In a recent paper, B. Green \cite{green} proved the following beautiful
result, which improves upon Bourgain's result above, and is currently
the best that is known on this problem: 

\begin{theorem} Suppose $A,B$ are subsets of ${\mathbb Z}/N{\mathbb Z}$
having cardinalities $\gamma N$ and $\delta N$, respectively.  Then
there is an absolute constant $c > 0$ such that $A+B$ contains an 
arithmetic progression of length at least
$$
\exp( c ((\gamma \delta \log N)^{1/2} - \log\log N)).
$$
\end{theorem}

There are also several other papers which treat the question of long
arithmetic progressions in sumsets $A+A+\cdots + A$, such as 
\cite{freiman}, \cite{sarkozy}, \cite{sarkozy2}, \cite{sarkozy3},
\cite{lev}, and \cite{lev2}. 

In this paper we give a very simple, elementary proof of a 
result, which shows that sumsets $A+B$ have long arithemtic progressions
when $A, B \subseteq \{1,2,...,N\}$ 
have only $N^{1-\epsilon}$ elements (the length of the
longest progression will depend on $\epsilon$).  This result is stronger
than those given in the above theorems of Bourgain and Green when $A$
and $B$ have less than $N/\log^A N$ elements, for some sufficiently large
$A$; however, when $A$ and $B$ have more than this many elements, their
results give a much stronger conclusion.  
The author would like to emphasize 
that the result given below is certainly not an improvement over the results
of Bourgain and Green, both of which use much more sophisticated harmonic
analytic methods, but it does show that it is possible to prove long
arithmetic progressions in very thin sumsets.   
\bigskip

\begin{theorem} \label{main_theorem}  
Suppose that 
$$
A,B \subseteq \{1,2,...,N\} 
$$
and 
$$
|A|,\ |B| > 5 N^{1-(4(k-1))^{-1}}.
$$
Then, the sumset $A+B$ must contain a non-trivial
$k$-term arithmetic progression, which is a sequence of
integers $n, n+d, n+2d,...,n+(k-1)d$, where $d \neq 0$.
\end{theorem}

To compare this result with those of Bourgain and Green, we note that
when $|A|, |B| \gg N$, then Green's result gives that $A+B$ contains 
a progression of length $\exp(c (\log N)^{1/2})$, for some constant $c$, 
whereas the author's result below will only give $\Omega(\log N)$.  
So, in this range, both Green's and Bourgain's result is much stronger than
the author's; however, when $|A|, |B| \ll N/\sqrt{\log N}$, 
then Green's result does
not give a non-trivial bound on the length of the longest 
arithmetic progression in $A+B$, whereas the author's result above gives
that $A+B$ contains a progression of length $\Omega( (\log N)/\log\log N)$.
\bigskip

\section{Proof of Theorem \ref{main_theorem}}
\bigskip

It is obvious that $A+B$ contains a $k$-term arithemtic progression
if and only if $A+B$ contains a $k$-term arithmetic progression modulo 
$4x + 1$; and, we will prove the theorem by showing that there is such
a progression modulo $4x+1$.  First, we need the following lemma:

\begin{lemma}  There exists a subset $C$ of the residue classes modulo $4x+1$
such that the following all hold:

1.  $|C| > 5x^{1-(k-1)^{-1}}$; 

2.  If $c \in C$, then there exists $d \in C$ such that $d \equiv -c \pmod{4x+1}$; and,

3.  If $C+C$ contains a non-trival $k$-term arithmetic progression modulo 
$4x+1$, then so must $A+B$.
\end{lemma}

\noindent {\bf Proof of the Lemma.}  We first claim that there exists an integer
$j \in \{1,2,...,4x\}$ such that there are at least $5 N^{1-(2k)^{-1}}$ residue
classes in common between the sets $A+j$ and $B$:  For a randomly selected
$j \in \{0,1,2,...,4x\}$ with the uniform probability measure, the probability that
a particular $b \in B$ happens to lie in $A+j = \{a + j \pmod{4x+1}\ :\ a \in A\}$
is at least $|A|/(4x+1)$; so, the expected size of the intersection of residue classes
between $A+j$ and $B$ is $|A||B|/(4x+1) > 5 N^{1-(2(k-1))^{-1}}$.  Since the average
intersection is this big, there must exist a $j$ for which $A+j$ and $B$
have at least $5 N^{1-(2(k-1))^{-1}}$ classes in common.  Let $D$ be the intersection
of such residue classes between $A+j$ and $B$.  We note that $D+D$ is a 
subset of $(A+j) +B = \{a+b+j\ :\ a \in A, b \in B\}$; and, it is obvious that if
$D+D$ contains a $k$-term arithemtic progression, then so must $A+B$.

Next, we show that there exists an integer $\ell \in \{0,1,...,4x\}$ such that
$D+\ell$ and $-(D+\ell)$ have at least $5 N^{1-(k-1)^{-1}}$ elements in common:
There are exactly the same number of residue classes common to 
$D+\ell$ and $-(D+\ell)$ are there are residue classes common to 
$D$ and $-D - 2\ell$.  Now, given $d \in D$, and a randomly chosen 
$\ell \in \{0,1,...,4x\}$ with the uniform measure, the probability that 
$d$ is congruent to some member of $-D - 2\ell$ is $|D|/(4x+1)$; 
so, on average, $D$ and $-D-2\ell$ have $|D|^2/(4x+1)$ residue classes
in common.  We deduce that there exists an integer $\ell$ such that 
$-D-2\ell$ and $D$ have at least $|D|^2 / (4x+1) > 5 N^{1-(k-1)^{-1}}$ classes in
common; and therefore, $D+\ell$ and $-(D+\ell)$ have at least 
$5 N^{1-(k-1)^{-1}}$ residue classes in common.  Now, let $C$ be the 
set of progressions common to $D+\ell$ and $-(D+\ell)$.  Then, trivially
properties 1 and 2 claimed by the lemma are satisfied.  Property 3 follows
since the residue classes occupied by $C+C$ are a translate of those
occupied by $A+B$, and since translations preserve arithmetic progressions.
\ \ \ \ $\blacksquare$
\bigskip

Resuming the proof of our theorem, we will show that $C+C$ contains a
$k$-term arithmetic progression.  First, we note that 
$$
x_1,\ x_2 = x_1 + d,\ ...,x_k = x_1 + (k-1)d
$$ 
are part of a $k$-term arithemtic progression modulo $4x+1$ if and only if
the following congruences are all satisfied
\begin{eqnarray} \label{congruence_system}
x_1 + x_3 &\equiv& 2x_2 \pmod{4N+1} \nonumber \\
x_2 + x_4 &\equiv& 2x_3 \pmod{4N+1} \nonumber \\
x_3 + x_5 &\equiv& 2x_4 \pmod{4N+1} \nonumber \\
&\vdots& \nonumber \\
x_{k-2} + x_k &\equiv& 2x_{k-1} \pmod{4N+1}. 
\end{eqnarray}

Now, since $C = -C$, we have that our set $C+C$ is the same as $C-C$.  
So, we may express the numbers $x_1,...,x_k \in C+C$ as 
$$
x_i\ =\ y_i - z_i,\ \ y_i,z_i \in C.
$$
Thus, we may re-express the congruences in (\ref{congruence_system}) as
\begin{eqnarray} \label{simplified_system}
y_1 + y_3 - 2y_2\ &\equiv&\ z_1 + z_3 - 2z_2 \pmod{4N+1} \nonumber \\
y_2 + y_4 - 2y_3\ &\equiv&\ z_2 + z_4 - 2z_3 \pmod{4N+1} \nonumber \\
y_3 + y_5 - 2y_4\ &\equiv&\ z_3 + z_5 - 2z_4 \pmod{4N+1} \nonumber \\
&\vdots & \nonumber \\
y_{k-2} + y_k - 2y_{k-1}\ &\equiv&\ z_{k-2} + z_k - 2z_{k-1} \pmod{4N+1}.
\end{eqnarray}

We now show that this system has ``lots'' of solutions:  Consider the
set of vectors of the form
\begin{equation} \label{vector}
(y_1 + y_3 - 2y_2 \pmod{4N+1},\ ...,\ y_{k-2} + y_k - 2y_{k-1} \pmod{4N+1} ),
\end{equation}
where $y_1,...,y_k \in C$.  Clearly, there can be at most $(4N+1)^{k-2}$ such
vectors, since there are $k-2$ coordinates and since each coordinate can take
on one of at most $4N+1$ different values.  On the other hand, there are 
$|C|^k$ different choices for $y_1,...,y_k$.  Now, each time we have a pair
of sequences $y_1,...,y_k$ and $z_1,...,z_k$ whose corresponding vector in
(\ref{vector}) is the same, we get a solution to (\ref{simplified_system}).  
To get a lower bound on the number of solutions to this system, we let 
$\lambda(n_1,...,n_{k-2})$ denote the number $k$-tuples 
$(y_1,...,y_k)$, each $y_i \in C$, such that 
$$
y_1 + y_3 - 2y_2 \equiv n_1 \pmod{4N+1},\ ...,\ 
y_{k-2} + y_k - 2y_{k-1} \equiv n_{k-2} \pmod{4N+1}.
$$
Then, the number of solutions $(y_1,...,y_k,z_1,...,z_k)$ to 
(\ref{simplified_system}) is 
$$
S\ =\ \sum_{0 \leq n_1,...,n_{k-2} \leq 4N} \lambda^2(n_1,...,n_{k-2}).
$$
Since there are $|C|^k$ choices for $y_1,...,y_k$, we get that 
$$
\sum_{0 \leq n_1,...,n_{k-2} \leq 4N} \lambda(n_1,...,n_{k-2})\ =\ |C|^k.
$$
So, the smallest that $S$ could be, subject to this constraint, is if all
the $\lambda(n_1,...,n_{k-2})$ were equal.  This gives the lower bound
$$
S\ >\ (4N+1)^{k-2} \left ( {|C|^k \over (4N+1)^{k-2}} \right )^2
\ >\ {|C|^{2k} \over (4N+1)^{k-2}}.
$$

Now, in order for this to allow us to conclude that there is a non-trivial
$k$-term AP, we must show that $S$ exceeds the number of solutions
to (\ref{simplified_system}) that give trivial solutions:  A trivial solution
occurs when $x_1 = \cdots = x_k$, and there are $|C|$ different ways that
this can happen, since $x_1$ (and $x_2,...,x_k$) can be any of the $|C|$ 
elements of $C$.  For each of these $|C|$ trivial progressions, there can
be at most $|C|$ different ways of writing $x_i = y_i - z_i$.  So, the total
number of solutions to (\ref{simplified_system}) that can lead to a
trivial progression is at most 
$$
|C| \times |C|^k\ =\ |C|^{k+1};
$$
and, this will be less than  
$$
{|C|^{2k} \over (4N+1)^{k-2}}\ <\ S,
$$
provided
$$
|C|\ >\ (4N+1)^{1 - {1 \over k-1}},
$$
which we know is satisfied.  Thus, the total number of solutions to 
(\ref{simplified_system}) outnumbers those solutions that can lead to
trivial arithemtic progressions, and so we conclude that $C+C$,
and therefore $A+B$, must have a non-trivial $k$-term AP.\ \ \ \ $\blacksquare$

\end{document}